\documentclass[reqno]{amsart}

\usepackage{amsmath,amsfonts,amssymb,amscd,verbatim,delarray,verbatim,epigraph}

\usepackage{amsthm}

\usepackage{url}

\newcommand{\F}{{\mathbb F}}

\newcommand{\Q}{{\mathbb Q}}

\begin{document}

\newtheorem{theorem}{Theorem}

\newtheorem{corollary}[theorem]{Corollary}
\newtheorem{corol}[theorem]{Corollary}
\newtheorem{conj}[theorem]{Conjecture}
\newtheorem{prop}[theorem]{Proposition}

\theoremstyle{definition}
\newtheorem{definition}[theorem]{Definition}
\newtheorem{example}[theorem]{Example}

\newtheorem{remarks}[theorem]{Remarks}
\newtheorem{remark}[theorem]{Remark}

\newtheorem{question}[theorem]{Question}
\newtheorem{problem}[theorem]{Problem}

\newtheorem{quest}[theorem]{Question}
\newtheorem{questions}[theorem]{Questions}

\def\toeq{{\stackrel{\sim}{\longrightarrow}}}
\def\into{{\hookrightarrow}}


\def\alp{{\alpha}}  \def\bet{{\beta}} \def\gam{{\gamma}}
 \def\del{{\delta}}
\def\eps{{\varepsilon}}
\def\kap{{\kappa}}                   \def\Chi{\text{X}}
\def\lam{{\lambda}}
 \def\sig{{\sigma}}  \def\vphi{{\varphi}} \def\om{{\omega}}
\def\Gam{{\Gamma}}   \def\Del{{\Delta}}
\def\Sig{{\Sigma}}   \def\Om{{\Omega}}
\def\ups{{\upsilon}}


\def\F{{\mathbb{F}}}
\def\BF{{\mathbb{F}}}
\def\BN{{\mathbb{N}}}
\def\Q{{\mathbb{Q}}}
\def\Ql{{\overline{\Q }_{\ell }}}
\def\CC{{\mathbb{C}}}
\def\R{{\mathbb R}}
\def\V{{\mathbf V}}
\def\D{{\mathbf D}}
\def\BZ{{\mathbb Z}}
\def\K{{\mathbf K}}
\def\XX{\mathbf{X}^*}
\def\xx{\mathbf{X}_*}

\def\AA{\Bbb A}
\def\BA{\mathbb A}
\def\HH{\mathbb H}
\def\PP{\Bbb P}

\def\Gm{{{\mathbb G}_{\textrm{m}}}}
\def\Gmk{{{\mathbb G}_{\textrm m,k}}}
\def\GmL{{\mathbb G_{{\textrm m},L}}}
\def\Ga{{{\mathbb G}_a}}

\def\Fb{{\overline{\F }}}
\def\Kb{{\overline K}}
\def\Yb{{\overline Y}}
\def\Xb{{\overline X}}
\def\Tb{{\overline T}}
\def\Bb{{\overline B}}
\def\Gb{{\bar{G}}}
\def\Ub{{\overline U}}
\def\Vb{{\overline V}}
\def\Hb{{\bar{H}}}
\def\kb{{\bar{k}}}

\def\Th{{\hat T}}
\def\Bh{{\hat B}}
\def\Gh{{\hat G}}


\def\cC{{\mathcal C}}
\def\cU{{\mathcal U}}
\def\cP{{\mathcal P}}
\def\cV{{\mathcal V}}
\def\cS{{\mathcal S}}

\def\CG{\mathcal{G}}

\def\cF{{\mathcal {F}}}

\def\Xt{{\widetilde X}}
\def\Gt{{\widetilde G}}


\def\hh{{\mathfrak h}}
\def\lie{\mathfrak a}

\def\XX{\mathfrak X}
\def\RR{\mathfrak R}
\def\NN{\mathfrak N}

\def\minus{^{-1}}

\def\GL{\textrm{GL}}            \def\Stab{\textrm{Stab}}
\def\Gal{\textrm{Gal}}          \def\Aut{\textrm{Aut\,}}
\def\Lie{\textrm{Lie\,}}        \def\Ext{\textrm{Ext}}
\def\PSL{\textrm{PSL}}          \def\SL{\textrm{SL}} \def\SU{\textrm{SU}}
\def\loc{\textrm{loc}}
\def\coker{\textrm{coker\,}}    \def\Hom{\textrm{Hom}}
\def\im{\textrm{im\,}}           \def\int{\textrm{int}}
\def\inv{\textrm{inv}}           \def\can{\textrm{can}}
\def\id{\textrm{id}}              \def\Char{\textrm{char}}
\def\Cl{\textrm{Cl}}
\def\Sz{\textrm{Sz}}
\def\ad{\textrm{ad\,}}
\def\SU{\textrm{SU}}
\def\Sp{\textrm{Sp}}
\def\PSL{\textrm{PSL}}
\def\PSU{\textrm{PSU}}
\def\rk{\textrm{rk}}
\def\PGL{\textrm{PGL}}
\def\Ker{\textrm{Ker}}
\def\Ob{\textrm{Ob}}
\def\Var{\textrm{Var}}
\def\poSet{\textrm{poSet}}
\def\Al{\textrm{Al}}
\def\Int{\textrm{Int}}
\def\Smg{\textrm{Smg}}
\def\ISmg{\textrm{ISmg}}
\def\Ass{\textrm{Ass}}
\def\Grp{\textrm{Grp}}
\def\Com{\textrm{Com}}
\def\rank{\textrm{rank}}

\def\char{\textrm{char}}

\def\wid{\textrm{wd}}

\newcommand{\Or}{\operatorname{O}}

\def\tors{_\def{\textrm{tors}}}      \def\tor{^{\textrm{tor}}}
\def\red{^{\textrm{red}}}         \def\nt{^{\textrm{ssu}}}

\def\sss{^{\textrm{ss}}}          \def\uu{^{\textrm{u}}}
\def\mm{^{\textrm{m}}}
\def\tm{^\times}                  \def\mult{^{\textrm{mult}}}

\def\uss{^{\textrm{ssu}}}         \def\ssu{^{\textrm{ssu}}}
\def\comp{_{\textrm{c}}}
\def\ab{_{\textrm{ab}}}

\def\et{_{\textrm{\'et}}}
\def\nr{_{\textrm{nr}}}

\def\nil{_{\textrm{nil}}}
\def\sol{_{\textrm{sol}}}
\def\End{\textrm{End\,}}

\def\til{\;\widetilde{}\;}

\def\min{{}^{-1}}

\def\AGL{{\mathbb G\mathbb L}}
\def\ASL{{\mathbb S\mathbb L}}
\def\ASU{{\mathbb S\mathbb U}}
\def\AU{{\mathbb U}}


\title{On first order rigidity for linear groups}

\author[Plotkin]{Eugene Plotkin}


\address{Plotkin: Department of Mathematics,
Bar-Ilan University, 5290002 Ramat Gan, ISRAEL}
\email{plotkin@macs.biu.ac.il}

\begin{abstract}

The paper is a short survey of recent developments in the area of
first order descriptions of linear groups. It is aimed to illuminate the known results and  to pose
the new problems  relevant to logical characterizations of Chevalley groups and Kac--Moody groups.
\end{abstract}

\maketitle

\epigraph{{\it To Prof. Jun Morita}\\
\smallskip
{\it with great respect and admiration}}

\bigskip



\section{Introduction}\label{Equiv}

Questions we are going to illuminate in this paper are concentrated around the interaction between algebra, logic, model theory and geometry.

The main question behind further considerations is as follows. Suppose we have two algebras equipped with a sort of logical description.

\begin{problem} When the coincidence of logical descriptions  provides an isomorphism between algebras in question?
\end{problem}

With this aim we consider different kinds of logical equivalences between algebras. Some of the notions we are dealing with are not formally defined in the text.
For  precise definitions and references  use \cite{Halmos}, \cite{MR}, \cite{Plotkin_UA-AL-Datab},  \cite{Plotkin_Haz}, \cite{Pl-St}, \cite{PlAlPl}, \cite{PlPl}.

First, we make emphasis on elementary equivalence of groups. Importance of the elementary classification of algebraic structures goes back to the famous  works of A.Tarski and A.Malcev. The main problem is to figure out {\it what are the algebras elementarily equivalent to a given one.} We will describe the current state of art of the problem: {\it when elementary equivalence of groups implies their isomorphism}. Situation of such kind will be called, for short, {\it elementary rigidity}.

Our second aim is to describe the notions of isotypicity of algebras and logical equivalence of algebras. These notions are much less known than elementary equivalence. However, they can logically characterize algebras in a very rigid way and one can expect affirmative answers to most of the problems formulated.

We discuss these notions from the perspectives of Chevalley groups and some other linear groups, and Kac-Moody groups.
\section{Elementary equivalence of algebras}\label{Equiv}

\subsection{Definitions}

Given an algebra $H$, its {\it elementary theory} $Th(H)$ is the set of all sentences (closed formulas) valid on $H$.

\begin{definition}\label{elem}
Two algebras $H_1$ and $H_2$ are said to be {\it elementarily equivalent} if their elementary theories coincide.
\end{definition}

Very often we fix a class of algebras $\mathcal C$ and ask what are the algebras elementarily equivalent to a given algebra inside the class $\mathcal C$. So, the rigidity question with respect to elementary equivalence looks as follows.


 \begin{problem}
 Let a class of algebras $\mathcal C$ and an algebra $H\in \mathcal C$ be given. Suppose that the elementary theories of algebras $H$ and $A\in \mathcal C$ coincide. Are they  {\it elementarily rigid}, that is,  are  $H$ and $A$  isomorphic?
 \end{problem}

For example, $\mathcal C$ can be the class of all groups, the class of finitely generated groups, the class of profinite groups, etc.

\begin{remark}
What we call elementary rigidity has different names. This notion appeared in the papers by A.Nies \cite{Nies}  under the name of quasi definability of groups. The corresponding name used in \cite{ALM} with respect to the class of finitely generated groups is first order rigidity. For some reasons which will be clear later on we use another term.
\end{remark}

In other words we ask for which algebras their logical characterization by means of the elementary theory is strong enough and defines the algebra in the unique, up to an isomorphism, way?

 We restrict our attention to the case of groups.  Elementary rigidity of groups occurs not very often. Usually various extra conditions are needed. Consider  examples of elementary rigidity for linear groups. First of all, a group which is elementarily equivalent to a finitely generated linear group is a residually finite linear group \cite{Mal1}. The incomplete list of known rigidity cases is given in the following theorem.

 \subsection{Chevalley groups}

 \begin{theorem}\label{mthm}
Historically, the first result was obtained by A.Malcev:
 \begin{itemize}
 \item  If two linear groups $GL_n(K)$ and  $GL_m(F)$, where $K$ and $F$ are    fields, are elementarily equivalent, then $n=m$ and the fields $K$ and $F$ are elementarily equivalent, see  \cite{Mal2}.
  \item This result was generalized to the wide class of Chevalley groups. Let $G_1=G_\pi(\Phi,R)$ and $G_2=G_\mu(\Psi,S)$ be two elementarily equivalent Chevalley groups. Here $\Phi$, $\Psi$ denote the root systems of rank $\geqslant 1$, $R$ and $S$ are local rings, and $\pi$, $\mu$ are weight lattices. Then root systems and weight lattices of $G_1$ and $G_2$ coincide, while the rings are elementarily equivalent. In other words Chevalley groups over local rings are elementarily rigid in the class of such groups modulo elementary equivalence of  ground rings \cite{Bu}.
      \item Let $G_\pi(\Phi,K)$ be a simple Chevalley group over the algebraically closed field $K$. Then  $G_\pi(\Phi,K)$ is elementarily rigid in the class of all groups (cardinality is fixed). This result can be deduced from \cite{Zilber}. In fact, this is true for a much wider class of algebraic groups over algebraically closed fields and, modulo elementary equivalence of fields,  over arbitrary fields \cite{Zilber}.
   \item Any  irreducible non-uniform higher-rank characteristic zero
arithmetic lattice is elementarily rigid in the class of all finitely generated groups, see \cite{ALM}. In particular, $SL_n(\mathbb Z)$, $n>2$ is elementarily rigid.
 \item Recently, the results of \cite{ALM} have been extended to a much wider class of lattices, see \cite{AM}.
   \item Let $\mathcal O$ be the ring of integers of a number field, and let $n\geqslant 3$. Then every  group $G$ which is elementarily equivalent to $SL_n(\mathcal O)$ is isomorphic to $SL_n(\mathcal R)$, where the rings $\mathcal O$ and $\mathcal R$ are elementarily equivalent. In other words $SL_n(\mathcal O)$ is elementarily rigid in the class of all groups modulo elementary equivalence of rings. The similar results are valid with respect to $GL_n(\mathcal O)$ and  to the triangular group   $T_n(\mathcal O)$ \cite{SM}. These results intersect in part with the previous items, since  the ring $R=\mathbb Z$ is elementarily rigid in the class of all finitely generated rings \cite{Nies}, and thus  $SL_n(\mathbb Z)$ is elementarily rigid in the class of all finitely generated groups.
         \item For the case of arbitrary Chevalley groups the results similar to  above cited are obtained in \cite{ST} by different machinery for a wide class of ground rings. Suppose  the Chevalley group $G=G(\Phi,R)$ of rank $\geqslant 2$ over the ring $R$ is given. Suppose that the ring $R$ is elementarily rigid in the class $\mathcal C$ of rings.  Then $G=G(\Phi,R)$ is elementarily rigid in the corresponding class $\mathcal C_1$ of groups  if $R$ is a field, $R$ is a local ring and $G$ is simply connected, $R$ is a Dedekind ring of  arithmetic type, that is the ring of $S$-integers of a number field, $R$ is Dedekind ring with at least 4 units and $G$ is adjoint. In particular, if  a ring of such kind is finitely generated then it gives rise to elementary rigidity of $G=G(\Phi,R)$ in the class of all finitely generated groups (see \cite{AKNS}).
              If $R$ of such kind is not elementarily rigid then  $G=G(\Phi,R)$ is elementarily rigid in the class of all groups modulo elementary equivalence of rings.
             \item The Chevalley group $G=G(\Phi,R)$ of rank $\geqslant 2$  is elementarily rigid in the class of all finitely generated groups, if $R$ is a ring of one-variable polynomials over the finite field, i.e., $R=F_q[x]$, $char F_q\neq 2$, see \cite{ST}, \cite{AKNS} and  \cite{CKPV}, where the bounded generation of such $G(\Phi,R)$ in elementary generators is proven.
                 \item The Chevalley group $G=G(\Phi,R)$   is elementarily rigid in the class of all finitely generated groups, if $R$ is a ring of Laurent polynomials over  the finite field, i.e., $R=F_q[x, x^{-1}]$, $char F_q\neq 2$. The proof also relies on \cite{ST}, \cite{AKNS} and  \cite{CKPV}.
                     \item  The Chevalley group $G=G(\Phi,R)$  is elementarily rigid in the class of all finitely generated groups, if $R$  is a finitely generated ring of $S$-integers in a global function field of
positive characteristic which has infinitely many units and satisfy some additional condition on $S$, see \cite{ST}, \cite{AKNS},  \cite{CKPV} and \cite{Queen}.
 \end{itemize}
  \end{theorem}

 For rank 1 simple linear groups the situation is quite different. For $R=\mathbb Z$  and $G=PSL_2(\mathbb Z)$ there is the following result of Z.Sela, cf., \cite{Se1}, see \cite{Kazach}.

  \begin{theorem}\label{sela_not_published}
  A finitely generated group $G$ is elementarily equivalent to $PSL(2, \mathbb Z)$
if and only if $G$ is a hyperbolic tower (over $PSL(2, \mathbb Z)$).
 \end{theorem}

 Let us make some comments regarding Theorem \ref{sela_not_published}. It was A.Tarski who asked whether one can distinguish between finitely generated free groups by means of their elementary theories. This formidable  problem has been solved in affirmative, that is all free groups have one and the same elementary theory. Moreover, all finitely generated groups elementarily equivalent to a given non-abelian free group have been explicitly described,  see \cite{KM}, \cite{Se}.
Paper \cite{Se1} extends this line.

 The result for $SL_2(\mathbb Z)$ can be deduced from Theorem \ref{sela_not_published}, see \cite{Kazach}.

  \begin{theorem}\label{sela_not_published1}
 A finitely generated group $G$ is elementarily equivalent to $SL_2(\mathbb Z)$ if and only if $G$ is the central extension of a hyperbolic tower over $PSL_2(\mathbb Z)$  by $\mathbb Z_2$  with the cocycle  $f: PSL_2(\mathbb Z) \times  PSL_2(\mathbb Z)\to \mathbb Z_2$, where $f(x,x)=1$ for all $x\in  PSL_2(\mathbb Z)$ of order 2, and $f(x,y)=0$, otherwise.
 \end{theorem}

 However, the situation in rank one Chevalley groups over the rings with infinitely many invertible elements looks different. Bounded generation of $SL_2(R)$ in elementary generators, where $R$ is a ring  of $S$-integers of a number field with infinitely
many units  is proved in \cite{MRS}. The same fact is true if $R$  is a ring of $S$-integers in a global function field of
positive characteristic which has infinitely many units and satisfies some additional condition on $S$, see \cite{Queen}. Hence, in these cases one can think about elementary rigidity.

Before going over the situation for Kac-Moody groups we shall cite an important model-theoretic result by F.Point on collaboration between Chevalley-Demazure functor and the operation of taking ultraproducts. Let $F$ be a non-principal ultrafilter on the set $I$ and let $G_\Phi(,)$ be a simple Chevalley-Demazure group scheme over $\mathbb Z$
defined by a root system $\Phi$. Let $K_i$, $i\in I$ be a collection of fields and  $G_i=G_\Phi(K_i)$ be the corresponding set of Chevalley groups. Denote by $\prod_FG_i$ the ultraproduct of the groups $G_i$ with respect to ultrafilter $F$. Theorem of F.Point (see \cite{Point}) basically says that the functors $G_\Phi(,)$ and $\prod_F$ commute. Namely,

\begin{theorem}\label{point}
$G_\Phi(\prod_FK_i)=\prod_FG_{\Phi}(K_i)$.
 Moreover the ultraproducts of the
unipotent, diagonal and monomial subgroups of $G_\Phi(K_i)$  are isomorphic to
the corresponding subgroups of $G_\Phi(\prod_FK_i)$.
\end{theorem}

Note that Theorem \ref{point} remains true (modulo minor conditions on the fields) also in the case of twisted Chevalley groups. From the perspectives of Kac-Moody groups and elementary rigidity of Chevalley groups, the following question is of great interest.

\begin{question} Does the statement of Theorem \ref{point} remain true for
\begin{itemize}
\item
local rings,
\item Dedekind rings of arithmetic type,
\item arbitrary rings.
\item Elementary subgroups of the Chevalley groups and rings as above.
\end{itemize}
\end{question}

 \subsection{Kac--Moody groups}
Given a generalized Cartan matrix $A$ and a field $k$ (or a ring
$R$), the value $G_A(k)$ of the Tits functor $G_A\colon \mathbb
Z$-$Alg\to Grp$ defines a minimal Kac--Moody group over $k$, see
\cite{Ti} (cf. \cite{Ca}, \cite{MR}). One can view this functor as a
generalization of the Chevalley--Demazure group scheme. We assume
that $A$ is indecomposable. As a rule, we assume that the functor
$G_A$ is simply connected. However, speaking about isomorphisms of affine Kac-Moody groups over fields and Chevalley groups over rings, we
will freely, often without special notice, use the common language
abuse, assuming that we go over to its subquotient, taking the
derived subgroup (resp. subalgebra) and factoring out the centre, if
necessary.

If $A$ is a definite matrix, the group $G_A(k)$ is a Chevalley group
$G_\Phi(k)$ where $\Phi$ is the root system corresponding to $A$.
These groups were considered  in the previous
section. If $A$ is of affine type, $G_A(k)$ is isomorphic to the Chevalley
group $G_\Phi(k[t,t^{-1}])$ where $k[t,t^{-1}]$ is the ring of
Laurent polynomials. The general case of Kac-Moody groups is covered by $A$ of indefinite type.
The first question is related to Theorem \ref{point}.

Let $F$ be a non-principal ultrafilter on the set $I$,  let $K_i$, $i\in I$ be a collection of fields and let the Tits functor $G_A\colon \mathbb
Z$-$Alg\to Grp$ defines a minimal Kac--Moody group over $K_i$. Denote by $\prod_FG_i$ the ultraproduct of the groups $G_i$ with respect to ultrafilter $F$.

\begin{question} Let $A$ be of affine type. Is it true that the Tits functor and ultraproducts commute, that is the formula
$
G_A(\prod_FK_i)=\prod_FG_A(K_i)
$
holds. Is it true that the same property is satisfied for arbitrary $A$. Consider separately the hyperbolic case.
\end{question}

The penultimate item of Theorem \ref{mthm} implies that

\begin{theorem} Let $G_A(k)$ be an affine Kac-Moody group over a finite field $k$. Then $G_A(k)$ is elementarily rigid in the class of all finitely generated groups.
\end{theorem}

Let now $G_A(k)$ be a Kac--Moody group of indefinite type. B.~R\'emy \cite{Remy1} and P.-E.~Caprace--B.~R\'emy \cite{CR1}
showed that the minimal indefinite adjoint Kac--Moody groups
$G_A(\mathbb F_q)$ are simple provided $q>n>2 $ where $n$ is the size of $A$. These groups are also simple for some matrices $A$ if  $n=2$ and $q>3$.  J.~Morita and B.~R\'emy \cite{MoR} proved that in the case
where $k$ is the algebraic closure of $\mathbb F_q$ the groups $G_A(k)$ are simple.
P.-E.~Caprace and K.~Fujiwara \cite{CF} showed that over finite
fields these (infinite) simple groups have infinite commutator
width. It seems extremely unlikely that these groups are elementarily rigid.

Let $G_A(k)$ be an incomplete Kac--Moody group. There are several
ways to complete this group with respect to an appropriate topology.

Let  $A$ be of affine type, that is, $\overline G_A(k)$ is
a complete affine Kac--Moody group.
Then $\overline G_A(k)$ is isomorphic to a Chevalley group of the
form $G_\Phi(k((t)))$ where $k((t))$ is the field of formal Laurent
series over $k$. By Theorem \ref{mthm} this group is elementarily rigid modulo elementary equivalence of the ground field. If the field $k$ is finite then this group is elementarily rigid in the class of finitely generated groups.

\section{Isotypic equivalence of algebras}\label{Equiviso}

 The aim of this section is to introduce another logical invariant  which describes algebras more rigidly than elementary equivalence.
Elementary equivalence of algebras $H_1$ and $H_2$ assumes  coincidence of all first order sentences valid on $H_1$ and $H_2$. What we are going to introduce requires coincidence of all types valid on $H_1$ and $H_2$. We call such a situation {\it isotypicity} of algebras.  Before going over results in this direction, we need to make some preparations.

 \subsection{Basics of universal algebraic geometry}\label{basic}






Fix a variety of algebras $\Theta$. Let $W(X)$, $X=\{x_1, \ldots, x_n\}$ denote the finitely generated  free algebra in $\Theta$. By equations in $\Theta$ we mean expressions of the form $w\equiv w'$, where $w$, $w'$ are words in $W(X)$ for some $X$. This is our first syntactic object.  Next, let $ \tilde \Phi=(\Phi(X), X\in \Gamma)$ be the multi-sorted Halmos algebra of first order logical formulas based on atoms $w\equiv w'$, $w$, $w'$ in $W(X)$, see \cite{Pl-St}, \cite{Plotkin_AGinFOL}, \cite{PlPl}. There is a special procedure to construct such an algebraic object which plays  the same role with respect to First Order Logic as Boolean algebras do with respect to Propositional calculus. One can view elements of  $\tilde\Phi=(\Phi(X), X\in \Gamma)$ just as first order formulas over $w\equiv w'$.


Let $X=\{x_1, \ldots , x_n \}$ and let $H$ be an algebra in the variety $\Theta$. We have an affine space $H^n=H^X$ of points  $\mu : X \to H$. For every $\mu$ we have also the $n$-tuple $(a_1, \ldots , a_n) = \bar a\in H^n$ with $a_i = \mu(x_i)$. For the given $\Theta$ we have the homomorphism $\mu : W(X) \to H$ and, hence, the affine space $H^n$ is viewed as the set of homomorphisms
$Hom(W(X),H).$
The classical kernel $Ker(\mu)$ corresponds to each point $\mu : W(X) \to H$. This is exactly the set of equations for which the  point $\mu$ is a solution. Every point $\mu$ has also the logical kernel $LKer(\mu)$, see \cite{PlAlPl}, \cite{Plotkin_7_lec},  \cite{Plotkin_AGinFOL}. Logical kernel $LKer(\mu)$ consists of all formulas $u \in \Phi(X)$ valid on the point $\mu$. This is always an ultrafilter in $\Phi(X)$.


So we define syntactic and semantic areas where logic and geometry operate, respectively. Connect them by a sort of Galois correspondence.

 Let $T$ be a system of equations in $W(X)$. The set $A=T'$ in the affine space $Hom(W(X),H)$ consisting of all solutions of the system $T$ corresponds  to $T$. Sets of such kind are called {\it algebraic sets}. Vice versa, given a set $A$ of points in the affine space consider all equations $T=A'$ having $A$ as the set of solutions. Sets $T$ of such kind are called {\it closed congruences} over $W$.

 We can do the same correspondence with respect to arbitrary sets of formulas. Given a set $T$ of formulas in algebra of formulas (set of elements) $\Phi(X)$, consider the set $A=T^L$ in the affine space, such that every point of $A$ satisfies every formula of $\Phi$. Sets of such
kind are called {\it definable sets}. Points of $A$ are called solutions of the set of formulas $T$.  Conversely,  given a set $A$ of points in the affine space consider all formulas (elements) $T=A^L$ having $A$ as the set of solutions. Sets $T$ of such kind are {\it closed filters} in  $\Phi(X)$. Given arbitrary $T$ and $A$ we can make their Galois closures $T^{''}$ and $A^{''}$, and $T^{LL}$ and $A^{LL}$.


\subsection{Logical equivalence of algebras}

All algebraic sets constitute a category with special rational maps as morphisms \cite{PlPl}. The same is true with respect to definable sets \cite{PlPl}. So, we can formulate logical closeness of algebras geometrically.



\begin{definition}\label{def:lsim} We call algebras $H_1$ and $H_2$  {\it logically similar}, if the categories of definable sets $LG_\Theta(H_1)$ and $LG_\Theta(H_2)$ are isomorphic.
\end{definition}








\begin{definition}\label{def:LG}
Algebras  are called  logically equivalent, if for every $X $ and every set of formulas $T$ in $\Phi(X)$   the  equality $T^{LL}_{H_1}=T^{LL}_{H_2}$ holds .
\end{definition}

The set $T^{LL}_{H}$ is called the {\it logical radical} of $T$ with respect to $H$. So algebras $H_1$ and $H_2$ are {\it logically equivalent} if for every set of formulas $T$ the logical radicals with respect to $H_1$ and $H_2$ coincide.

First of all, it is easy to see that if algebras $H_1$ and $H_2$ are
logically equivalent then they are logically similar.  Now we want to understand what is the meaning of logical equivalence.

\begin{definition}\label{def:iso}
 Two algebras $H_1$ and $H_2$ are called LG-isotypic  if for every point $\mu:W(X)\to H_1$ there exists a point $\nu:W(X)\to H_2$ such that $LKer(\mu)=LKer(\nu)$ and, conversely, for every point $\nu:W(X)\to H_2$ there exists a point $\mu:W(X)\to H_1$ such that $LKer(\nu)=LKer(\mu)$.
\end{definition}

We can reformulate isotypicity of algebras in more standard logical notations.

\begin{definition}
Let $\mathcal{L}$ be a first-order language,  ${H}$ and ${G}$ be $\mathcal{L}$-algebras. Then ${H}$ and ${G}$ are isotypic,
if for any finite tuple $\bar{a}$ in ${H^n}$, there exists a tuple $\bar{b}$ in ${G^n}$ such that $tp^{H}(\bar{a})=tp^{G}(\bar{b})$ and vice versa.
\end{definition}

The meaning of the Definition \ref{def:iso} is the following. Two algebras are isotypic if the sets of realizable types over $H_1$ and $H_2$ coincide. So, by some abuse of language these algebras have the same logic of types. Some references for the notion of isotypic algebras are contained in \cite{Plotkin_M}, \cite{Plotkin_Gagta}, \cite{PlAlPl}, \cite{PlPl}, \cite{PZ}, \cite{Zhitom_types}. Note that the notion was introduced in \cite{PZ}, \cite{Plotkin_M} while \cite{PlPl} gives the most updated survey.

The principal property is as follows, see \cite{Zhitom_types}.
\begin{theorem}\label{thm:lgiso} 
Algebras $H_1$ and $H_2$ are logically equivalent if and only if they are isotypic.
\end{theorem}

\begin{definition}\label{def:lsim} We call the condition $\mathcal A$ rigid (or $\mathcal A$-rigid) in the class of algebras $\mathcal C$ if two algebras $H_1$ and $H_2$ from $\mathcal C$ subject to  $\mathcal A$ are isomorphic.
\end{definition}

Now we are in a position to study   rigidity of algebras with respect to isotypicity property.  It is easy to see that  
\begin{prop}
 If algebras $H_1$ and $H_2$ are
logically equivalent then they are elementarily equivalent.
\end{prop}

\subsection{Isotypic algebras}

 It is clear, that since isotypicity is stronger than elementary equivalence, this phenomenon can occur quite often. Let us state this problem explicitly.

 \begin{problem}
Let a class of algebras $\mathcal C$ and an algebra $H\in \mathcal C$ be given. Suppose that algebras $H\in \mathcal C$ and $A\in \mathcal C$ are isotypic. Are they  {\it isotypically rigid}, that is are $H$ and $A$  isomorphic?
 \end{problem}

\begin{remark}
In many papers  isotypically rigid algebras are called {\it logically separable} \cite{PlPl}, \cite{Plotkin_Gagta}, or {\it type definable} \cite{MRoman}. 
\end{remark}

The following principal  problem was stated in \cite{PlPl}, \cite{MRoman} and is widely open.

\begin{problem}[Rigidity problem]\label{princ}
 Is it true that every two isotypic finitely generated groups are
isomorphic?
\end{problem}

Meantime, Problem \ref{princ} is answered in affirmative for many groups. Some of the cases are collected in Theorem \ref{isor} and the consequent Corollary.

 \begin{theorem}\label{isor}
 The following cases of isotypically rigid groups are known:

 \begin{itemize}
 \item Every finitely generated co-Hopfian group is isotypically rigid in the class of all groups, see \cite{Zhitom_types}, \cite{Sklinos}.

 \item Every finitely presented Hopfian group is isotypically rigid in the class of all groups, see \cite{Sklinos}.

     \item Let $\Theta$ be a variety of groups. If a finitely generated free group in $\Theta$ is Hopfian then it is isotypically rigid in the class of all groups, see \cite{Zhitom_types}.

  \item  Finitely generated   metabelian  groups   are isotypically rigid in the class of all groups \cite{MRoman}.
 \item Finitely generated  virtually polycyclic groups are isotypically rigid in the class of all groups \cite{MRoman}.
 \item Finitely generated torsion free hyperbolic groups are isotypically rigid in the class of all groups \cite{Sklinos}.
 \item All surface groups, which are not  non-orientable surface
groups of genus 1,2 or 3   are isotypically rigid in the class of all groups \cite{MRoman}.

  \end{itemize}
  \end{theorem}

 \begin{corollary}\label{iso}
 Finitely generated absolutely free, free abelian, free nilpotent, free solvable groups are isotypically rigid.
 \end{corollary}

 \begin{conj}\label{con_lin_iso} Every finitely generated linear group is isotypically rigid. 
 \end{conj}

 \begin{conj}\label{con_lin_iso1} Let a Chevalley group $G_\Phi(R)$  over a ring $R$  be isotypic to a group $H$. Then $H$ is isomorphic to $G_\Phi(S)$ such that $R$ and $S$ are isotypic rings. 
 \end{conj}

\begin{problem}
What are the isotypicity classes of fields? When  two isotypic fields are isomorphic?
\end{problem}

\begin{remark}
In fact, using either logical equivalence of algebras, or what is the same, the isotypicity of algebras, we compare the possibilities of individual points in the affine space  to define the sets of formulas (in fact ultrafilters in $\Phi(X)$) which are valid in these points. Given a point $\mu$ in the affine space,  the collection of formulas valid in the point $\mu$ is {\it a type} of $\mu$. If these individual types are,  roughly speaking, the same for both algebras,  then these algebras are declared isotypic. Thus, for isotypic algebras  we compare  types of formulas realizable on these  algebras. Of course, this is significantly stronger than elementary equivalence, where the individuality of points disappeared and we compare only formulas valid in all points of the affine space.
\end{remark}

{\bf Acknowledgements.} The research  was supported by  ISF grant 1994/20  and the Emmy Noether Research Institute for Mathematics.




\end{document}